%% file: main.tex
\documentclass[11pt,reqno,a4paper]{article}
\usepackage[a4paper, total={6in, 8in}]{geometry}
\usepackage{geometry}
\usepackage{amsmath,amsthm,amstext,amscd,amssymb,euscript,url}
\usepackage{mathrsfs}
\usepackage{calrsfs}
\usepackage{epsfig}
\usepackage[inline]{enumitem}
\usepackage{todonotes}
\usepackage{mathtools}
\usepackage[absolute,overlay]{textpos}
\usepackage{graphicx}
\usepackage{subfig}
\usepackage{tikz}
\usetikzlibrary{matrix}
\usetikzlibrary{arrows,positioning}
\usepackage{color}

\usetikzlibrary{3d}
\usepackage{xifthen}
\usetikzlibrary{decorations.pathmorphing}

\input{mymacros}

\begin{document}

\title{{\bf Mosco convergence of independent particles and applications to particle systems with self-duality}}
\author{Mario Ayala \\
\small{School of Computation, Information and Technology}\\
\small{Chair for Analysis and Modelling}\\
\small{Technische Universität München}\\
{\small Boltzmannstraße 3}\\
{\small 85747 Garching }
\\
\small{Germany}
}
\maketitle

\begin{abstract}
We consider a sequence of Markov processes $\{X_t^n \mid n \in \N \}$ with Dirichlet forms converging in the Mosco sense of Kuwae and Shioya to the Dirichlet form associated with a Markov process $X_t$. Under this assumption, we demonstrate that for any natural number $k$, the sequence of Dirichlet forms corresponding to the Markov processes generated by $k$ independent copies of $\{X_t^n \mid n \in \N \}$ also converges. As expected, the limit of this convergence is the Dirichlet form associated with $k$ independent copies of the process $X_t$. We provide applications of this result in the context of \textbf{interacting} particle systems with Markov moment duality.
\end{abstract}

\section{Introduction}
Understanding the convergence of sequences of processes is fundamental for various applications, including in particular the micro to macro transition of statistical mechanics. This paper focuses on the Mosco convergence of Dirichlet forms in the setting of varying Hilbert spaces introduced by Kuwae and Shioya in \cite{kuwae2003convergence}, which, despite being weaker than weak convergence of stochastic processes, can play a crucial role in the convergence analysis of Markov processes. The applicability of Mosco convergence has already been shown in many works, including \cite{kolesnikov2006mosco,andres2010particle,10.1214/23-AAP2007,barret2014averaging}. In particular, we would like to mention the paper \cite{ayala2021condensation}, where the limiting variance of the density fluctuation field has been derived in a regime as complex as condensation of particles, for which weak convergence of the underlying Markov process was not directly accessible.

In this work, we consider a sequence of Markov processes $\lbrace X^n_t \mid n \in \mathbb{N}\rbrace$ with associated Dirichlet forms converging in the Mosco sense of \cite{kuwae2003convergence} to the Dirichlet form associated with a limiting Markov process $X_t$. Our main result extends this convergence to systems consisting of $k$ independent copies of these processes. Specifically, we demonstrate that if the Dirichlet forms of the original processes converge in the Mosco sense, then the Dirichlet forms of $k$ independent copies of these processes also converge to the Dirichlet form of $k$ independent copies of the limiting process $X_t$. This result is, of course, natural; however, to the best of our knowledge, there is no proof in the literature. This foundational result becomes particularly useful for interacting particle systems that exhibit the property of self-duality.

Self-duality provides a powerful tool for analyzing interacting systems of particles by relating the dynamics of a process with a simpler version of itself. This property simplifies the study of long-term behavior and correlations. For example, duality has been proven to be useful in a variety of contexts \cite{franceschini2024orthogonal,casini2024density,kuan2023dualities,wagner2024orthogonal,ayala2021higher}, among many others. By utilizing self-duality, we can reduce problems involving an infinite number of particles to the analysis of only a finite number of them. Hence the need for the derivation of convergence results for $k$-particles dynamics. As a byproduct, we show that our convergence result has significant implications for particle systems with Markov moment duality. In particular, thanks to Dirichlet domination, we are able to derive the Mosco convergence of the Dirichlet forms associated with a system of $k$ interacting particles with dynamics given by the symmetric inclusion process ($\SIP$) to the corresponding Brownian Dirichlet forms. Once equipped with this result, we derive weak scaling limits for the $k$-point correlation functions (the so-called $v$-functions) for the system consisting of an infinite number of $\SIP$ particles  started out of stationarity.

The rest of the paper is organized as follows: Section 2 introduces the necessary preliminaries on Hilbert space convergence and Mosco convergence. Section 3 presents the main results and assumptions. Section 4 provides detailed proofs of the main results. Finally, in Section 5, we present the derivation of scaling limits for the $k$-point correlation functions for the $\SIP$ diffusively rescaled.

\section{Preliminaries}\label{Preeliminaries}

\subsection{Convergence of Hilbert spaces}
The convergence of Dirichlet forms we are interested in takes place in the setting of convergence of a sequence of Hilbert spaces. We recall this notion of convergence introduced in \cite{kuwae2003convergence}.

\bd[Convergence of Hilbert spaces]\label{HilConv}
A sequence of Hilbert spaces $\{ H_n \}_{n \geq 0}$ converges to a Hilbert space $H$ if there exist a dense subset $\Cfrak \subseteq H$ and a family of linear maps $\left\{ \Phi_n \colon \Cfrak \to H_n \right\}_n$ such that:
\be\label{HilCond}
\lim_{n \to \infty} \| \Phi_n f \|_{H_n} = \| f \|_{H}, \qquad \text{  for all } f \in \Cfrak. 
\ee
\ed
\noindent
It is also necessary to introduce the concepts of strong and weak convergence of vectors living on a convergent sequence of Hilbert spaces. Hence in Definitions \ref{strongcon}, \ref{weakcon} and \ref{MoscoDef} below we assume that the spaces $\{ H_n \}_{n \geq 0}$ converge to the space $H$, in the sense we have just defined,  with the dense set  $\Cfrak \subset H$ and the sequence of operators $\{ \Phi_n \colon \Cfrak  \to H_n \}_n$ witnessing the convergence.

\bd[Strong convergence on Hilbert spaces]\label{strongcon}
A sequence of vectors $\{  f_n \}$ with $f_n$ in $H_n$, is said to \textbf{strongly converge} to a vector $f \in H$ if there exists a sequence $\{ \tilde{f}_M \} \in \Cfrak$ such that
\be\label{DefStrongcond1}
\lim_{M\to \infty} \| \tilde{f}_M -f \|_{H} = 0 
\ee
and
\be\label{DefStrongcond2}
\lim_{M \to \infty} \limsup_{n \to \infty} \| \Phi_n \tilde{f}_M -f_n \|_{H_n} = 0. 
\ee
\ed
\bd[Weak convergence on Hilbert spaces]\label{weakcon}
A sequence of vectors $\{  f_n \}$ with $f_n \in H_n$, is said to \textbf{weakly converge}  to a vector $f$ in a  Hilbert space $H$ if
\be
\lim_{n\to \infty} \left \langle f_n, g_n \right \rangle_{H_n} = \ \left \langle f, g \right \rangle_{H}, \nn 
\ee
for every sequence $\{g_n \}$ strongly convergent to $g \in H$.
\ed
\br\label{StrongConvPhin}
Notice that, as expected, strong convergence implies weak convergence, and, for any $f \in \Cfrak$, the sequence $\Phi_n f $ strongly-converges to $f$.
\er

\noindent
Given these notions of convergence, we can also introduce related notions of convergence for operators. More precisely, if we denote by $L(H)$ the set of all bounded linear operators in $H$, we  have the following definition.

\bd[Convergence of bounded operators on Hilbert spaces]\label{opercon}
A sequence of bounded operators $\{  T_n \}$ with $T_n \in L(H_n)$, is said to strongly (resp. weakly ) converge to an operator  $T$ in $L(H)$ if for every strongly (resp. weakly) convergent sequence $\{  f_n \}$, $f_n \in H_n$ to $f \in H$ we have that the sequence $\{ T_n f_n \}$ strongly (resp. weakly ) converges to $T f$.
\ed

\noindent
We are now ready to introduce Mosco convergence.

\subsection{Definition of Mosco convergence}
In this section we assume the Hilbert convergence of a sequence of Hilbert spaces $\{ H_n \}_n$ to a space $H$.

\bd[Mosco convergence]\label{MoscoDef}
A sequence of Dirichlet forms $\{ (\caE_n, D(\caE_n))\}_n $, defined on Hilbert spaces $H_n$, Mosco converges to a Dirichlet form $(\caE, D(\caE)) $, defined in some Hilbert space $H$, if:
\begin{description}
\item[Mosco I.] For every sequence of $f_n \in H_n$ weakly-converging  to $f$ in $H$
\be\label{mosco1}
\caE ( f ) \leq \liminf_{n \to \infty} \caE_n ( f_n ).
\ee
\item[Mosco II.] For every $f \in H$, there exists a sequence $ f_n \in H_n$ strongly-converging  to $f$ in $H$, such that
\be\label{mosco2}
\caE ( f) = \lim_{n \to \infty} \caE_n ( f_n ).
\ee
\end{description}
\ed

The following theorem from \cite{kuwae2003convergence}, which relates Mosco convergence with convergence of semigroups and resolvents, is a powerful application, and one of the main ingredients of the applications of our main result.
\bt\label{MKS}
Let $\{ (\caE_n, D(\caE_n))\}_n $ be a sequence of Dirichlet forms on Hilbert spaces $H_n$ and let  $(\caE, D(\caE)) $ be a Dirichlet form  in some Hilbert space $H$. The following statements are equivalent:
\ben
\item $\{ (\caE_n, D(\caE_n))\}_n $ Mosco-converges to $\{ (\caE, D(\caE))\} $.
\item The associated sequence of semigroups $\{ T_{n} (t) \}_n $ strongly-converges to the semigroup $ T(t)$ for every $t >0$.
\een
\et

\section{Main result}

\subsection{Assumptions}
Let $\left\{ \, X_t^{(n)} \mid n \in \N, \, t \geq 0 \, \right\} $ be a sequence of reversible Markov processes with state space $\Omega_n$ and reversible measures $\left\{ \, \mu_n \mid n \in \N \, \right\} $. Denote by $\left\{ \, (\caE_n, D(\caE_n)) \mid n \in \N \, \right\}$ the Dirichlet forms associated to $\left\{ \, X_t^{(n)} \mid n \in \N \, \right\} $ and defined on the Hilbert spaces $\left\{ \, L^2(\Omega_n, \mu_n) \mid n \in \N \, \right\}$. We also consider a  reversible Markov process $\left\{ \, X_t \mid t \geq 0 \, \right\} $ with reversible measure $\mu$ on a state space $\Omega$, and denote by $\left\{ \, \caE, D(\caE) \, \right\}$ its Dirichlet form.

\begin{assumption}\label{MainAssumption}
Let us assume the following.
\begin{enumerate}
    \item There exist a dense linear subspace $\Cfrak \subseteq D(\caE) \subseteq L^2(\Omega, \mu)$ and a family of linear maps 
    \[
    \left\{ \, \Phi_n \colon \Cfrak \to L^2(\Omega_n, \mu_n)\mid n \in \N \, \right\}
    \]
    such that the sequence of spaces $\left\{ \, L^2(\Omega_n, \mu_n) \mid n \in \N \, \right\} $ converges in the sense of \mbox{Definition~\ref{HilConv}} to the space $L^2(\Omega, \mu)$ with the couple $\left\{ \, \Cfrak, \{ \Phi_n\}_{n \in \N}\, \right\}$ witnessing the Hilbert convergence.
    \item The sequence of Dirichlet forms  $\left\{ \, (\caE_n, D(\caE_n)) \mid n \in \N \, \right\}$ converges in the sense of \mbox{Definition~\ref{MoscoDef}} to the Dirichlet form $\left\{ \, \caE, D(\caE) \, \right\}$. This convergence takes place under the Hilbert convergence 
    \[
    \left\{ \, L^2(\Omega_n, \mu_n) \mid n \in \N \, \right\} \rightarrow L^2(\Omega, \mu).
    \]
    \item There exists a core $\Cfrak_c \subseteq \Cfrak \subseteq D(\caE)$.
\end{enumerate}
\end{assumption}

\subsection{The $k$ independent motion}
Let $k \in \N$ and consider the sequence of processes  $\left\{ \, X_t^{(n,k)} \mid n \in \N, \, t \geq 0 \, \right\} $ consisting of $k$ independent copies of the processes $\left\{ \, X_t^{(n)} \mid n \in \N, \, t \geq 0 \, \right\} $. The state space of these processes is $\Omega_n^{(k)}:= \bigotimes_{j=1}^k\Omega_n$. In the same way we define the process $\left\{ \, X_t^{(k)} \mid t \geq 0 \, \right\} $ and its state space $\Omega^{(k)}$. Moreover, it is clear that each of the processes $\left\{ \, X_t^{(n,k)} \mid  t \geq 0 \, \right\}$ is reversible with respect to the product measure $\mu_n^{(k)}:=\otimes_{j=1}^k \mu_n$.

\br\label{IdentRemark} 
It is known, see for example 2.6.11 in \cite{kadison1986fundamentals} for $k=2$, that we have the identification
\be\label{TensProdHilIndent}
L^2(\Omega^{(k)}, \mu^{(k)}) \cong \bigotimes_{j=1}^k  L^2(\Omega, \mu).
\ee
In order to avoid issues of overloading notation as a consequence of specifying the operator witnessing \eqref{TensProdHilIndent} we choose once and for all an orthonormal basis $\left\{  e_i \mid i \in \N \right\}$ in $L^2(\Omega, \mu)$.
\er 
Notice that there is a natural inner product on $\bigotimes_{j=1}^k  L^2(\Omega, \mu)$ given by 
\be 
\langle \bif, \mathbf{g} \rangle = \prod_{j=1}^k \langle f_i, g_i \rangle
\ee 
for all $\mathbf{f}=(f_1, f_2, \ldots, f_k)$, $\mathbf{g}=(g_1, g_2, \ldots, g_k)$ with each $f_j, g_j \in  L^2(\Omega, \mu)$. Moreover, we also have
\be\label{DefHnk}
\norm{\mathbf{f}}_{k} := \prod_{j=1}^k \norm{f_j}, 
\ee 
and as a consequence the Dirichlet form associated to $\left\{ \, X_t^{(k)} \mid t \geq 0 \, \right\} $ is given by
\be\label{KDirFormDef}
\caE^{(k)}(\mathbf{f})= \sum_{j=1}^k \norm{\bif_{-j}}^2 \caE(f_j), 
\ee 
for all $\mathbf{f}=(f_1, f_2, \ldots, f_k)$ with each $f_j \in  L^2(\Omega, \mu)$ and where $\mathbf{f}_{-j}=(f_1, \ldots, f_{j-1}, f_{j+1}, \ldots, f_k)$.
Using the $n$-version of \eqref{TensProdHilIndent} we can define the Dirichlet form associated to  $\left\{ \, X_t^{(n,k)} \mid  \, t \geq 0 \, \right\}$ as
\be\label{KnDirFormDef}
\caE_n^{(k)}(\mathbf{f})= \sum_{j=1}^k \norm{\bif_{-j}}^2 \caE_n(f_j).
\ee 

\subsubsection{Convergence of tensor product Hilbert spaces}
Let $k \in \N$, we now introduce the $k$-versions of $\Cfrak$ and $\Phi_n$. We start with the first one.
\be\label{DefCfrakk}
\Cfrak^{(k)} := \operatorname{Span}\left\{ \, f_1 \otimes f_2 \ldots \otimes f_k  \mid f_j \in \Cfrak  \text{ for all } j \in [k]    \, \right\}.
\ee 

The following proposition confirms that $\Cfrak^{(k)}$ is indeed dense in $L^2(\Omega^{(k)}, \mu^{(k)})$.

\bp\label{DensityofCfrakk} 
$\Cfrak^{(k)}$ is a dense subset of $L^2(\Omega^{(k)}, \mu^{(k)})$.
\ep 

\bpr
We do not include the proof of this result since it follows a standard $\epsilon/4$ argument combining the density of $\Cfrak$ in $L^2(\Omega, \mu)$ and that of the orthonormal basis $\left\{  e_i \mid i \in \N \right\}$.
\epr 

We can now define the operators $\Phi_n^{(k)} \colon \Cfrak^{(k)} \to L^2(\Omega_n^{(k)}, \mu_n^{(k)})$. Take $\mathbf{f} \in \Cfrak^{(k)}$, i.e., there exist $M \in \N$ and $c_i \in \R$ such that 
\be\label{Expanbif} 
\mathbf{f} = \sum_{i=1}^M \alpha_i \otimes_{j=1}^k f_{(i,j)}, 
\ee 
for some $f_{(i,j)} \in \Cfrak$. Then we define the operator as 
\be\label{Defphink}
\Phi_n^{(k)} \mathbf{f}  := \otimes_{j=1}^k \Phi_n  \sum_{i=1}^M \alpha_i  f_{(i,j)}.
\ee 

\br
Notice that from the assumption that $\Cfrak$ is a vector space, and the multilinearily of the tensor product, we can simply substitute \eqref{Expanbif} by the simpler expansions
\be\label{SimplExpanbif} 
\bif = \otimes_{j=1}^k f_{(i,j)}. 
\ee 
We will use these simpler expansions whenever is convenient for notational purpuses.
\er

\bl\label{LemmaConvKHilb} 
For every $k \in \N$, let $\left\{ \, \Cfrak^{(k)}, \{ \Phi_n^{(k)}\}_{n \in \N}\, \right\}$ be defined as in \eqref{DefCfrakk} and \eqref{Defphink}. We then have
\be 
    \left\{ \, L^2(\Omega_n^{(k)}, \mu_n^{(k)}) \mid n \in \N \, \right\} \rightarrow L^2(\Omega^{(k)}, \mu^{(k)}), \nn 
\ee 
in the sense of \mbox{Definition~\ref{HilConv}} with the couple $\left\{ \, \Cfrak^{(k)}, \{ \Phi_n^{(k)}\}_{n \in \N}\, \right\}$ witnessing the convergence.
\el 

\bpr[Lemma \ref{LemmaConvKHilb}]
It is enough to verify \eqref{HilCond}. Let $\mathbf{f} = \sum_{i=1}^M \alpha_i \otimes_{j=1}^k f_{(i,j)} \in \Cfrak^{(k)}$, then we have
\begin{align}
    \lim_{n \to \infty} \norm{\Phi_n^{(k)} \mathbf{f}}_{n,k}^2 &=  \lim_{n \to \infty} \norm{\otimes_{j=1}^k \Phi_n  \sum_{i=1}^M \alpha_i  f_{(i,j)}}_{n,k}^2 =  \lim_{n \to \infty} \prod_{j=1}^k \norm{ \Phi_n  \sum_{i=1}^M \alpha_i  f_{(i,j)}}_{n,k}^2 \nn \\ 
    &=  \prod_{j=1}^k \norm{\sum_{i=1}^M \alpha_i  f_{(i,j)}}_{\infty,k}^2 =  \norm{ \otimes_{j=1}^k \sum_{i=1}^M \alpha_i  f_{(i,j)}}_{\infty,k}^2 = \norm{\mathbf{f}}_{\infty,k}^2, 
\end{align}
where in the second equality we used definition \eqref{DefHnk}.
\epr 

The following is the main result of our work.
\bt\label{MainThmkpart}
Under the Hilbert converges of Lemma \ref{LemmaConvKHilb}, the sequence of processes $\left\{ \, X_t^{(n,k)} \mid n \in \N, \, t \geq 0 \, \right\} $ given by the $k$-independent motion of $\left\{ \, X_t^{(n)} \mid n \in \N, \, t \geq 0 \, \right\} $ converges in the sense of Mosco convergence of Dirichlet forms to the Dirichlet forms associated to the process $\left\{ \, X_t^{(k)} \mid t \geq 0 \, \right\}$ given by the $k$-independent motion of the process $\left\{ \, X_t \mid t \geq 0 \, \right\}$.
\et

\section{Proofs of main results}

\subsection{Proof of condition Mosco I}
We start with the following auxiliary result.
\bp\label{PropweakConv1fromk} 
Let $\left\{ \,  \bif_n=\otimes_{j=1}^k f^{(j)} \mid n \in \N \, \right\}$ be a sequence of elements of the Hilbert space $L^2(\Omega_n^{(k)}, \mu_n^{(k)})$ converging weakly to a non-zero element $\bif=\otimes_{j=1}^k f^{(j)}$ of $ L^2(\Omega^{(k)}, \mu^{(k)})$. Then for every $j \in [k]$, we have that the sequence $\left\{ \,  f_n^{(j)} \mid n \in \N \, \right\}$ converges weakly to $f^{(j)}$.
\ep

\bpr[Proposition \ref{PropweakConv1fromk}] 
We only show the proof for the case $k=2$, since other cases are completely analogous. Let $\left\{ \,  f_n^{(1)} \otimes f_n^{(2)} \mid n \in \N \, \right\}$  converge weakly to $f^{(1)} \otimes f^{(2)}$, where both \mbox{$f^{(j)}\neq 0, j \in {1,2}$}. It is enough to show that $\left\{ \,  f_n^{(1)}  \mid n \in \N \, \right\}$  converges weakly to $f^{(1)}$. Let then $\left\{ \,  g_n^{(1)}  \mid n \in \N \, \right\}$ be a sequence converging strongly to some $g^{(1)} \in L^2(\Omega, \mu)$. That is, there exists a sequence $\tilde{g}_m^{(1)} \in \Cfrak$ such that\be\label{Strong1condgm}
\lim_{m \to \infty} \norm{ \tilde{g}_m^{(1)}-g^{(1)}}_{\infty} = 0 
\ee
and
\be\label{Strong2condgm}
\lim_{m \to \infty} \limsup_{n \to \infty} \| \Phi_n \tilde{g}_m^{(1)}-g_n^{(1)}\|_{n} = 0. 
\ee
The idea is to construct a sequence $\bfg_n $ in $L^2(\Omega_n^{(2)}, \mu_n^{(2)})$ converging strongly to some $\bfg$ in $L^2(\Omega^{(2)}, \mu^{(2)})$ that allows us to recover the convergence
\be 
\lim_{n\to \infty} \left \langle f_n^{(1)}, g_n^{(1)} \right \rangle_{n} = \ \left \langle f^{(1)}, g^{(1)} \right \rangle_{\infty}. \nn
\ee 
The natural choices are
\be 
\bfg_n := g_n^{(1)} \otimes g_n^{(2)} := g_n^{(1)} \otimes \frac{f_n^{(2)}}{\norm{f_n^{(2)}}_n^2},
\ee 
and 
\be 
\bfg := g^{(1)} \otimes g^{(2)} := g^{(1)} \otimes \frac{f^{(2)}}{\norm{f^{(2)}}_\infty^2},
\ee 
whenever $f_n^{(2)}$ is not zero(which happens only finetely many times), otherwise $f_n^{(2)}$ is just replaced by $\Phi_n f^{(2)}$\\

We now show that indeed $\bfg_n$ converges strongly to $\bfg$ by checking that \eqref{DefStrongcond1} and \eqref{DefStrongcond2} are satisfied. 

In order to verify \eqref{DefStrongcond1} it is enough to notice that by the density of $\Cfrak$ there exists a sequence  $\tilde{g}_m^{(2)} \in \Cfrak$ converging in norm to $g^{(2)}$ and define $\tilde{\bfg}_m:=\tilde{g}_m^{(1)} \otimes \tilde{g}_m^{(2)}$. We then have:
\begin{align}
\norm{ \tilde{\bfg}_m-\bfg}_{\infty,2} &= \norm{ \tilde{g}_m^{(1)}-g^{(1)}}_{\infty} \norm{ \tilde{g}_m^{(2)}-g^{(2)}}_{\infty}
\end{align}
which by \eqref{Strong1condgm}, and the norm convergence $\tilde{g}_m^{(2)} \to g^{(2)}$, implies \eqref{DefStrongcond1}.

For condition \eqref{DefStrongcond2}, we simply estimate
\begin{align}
    \| \Phi_n^{(2)} \tilde{\bfg}_m-\bfg_n \|_{n, 2}^2 &= \| \Phi_n \tilde{g}_m^{(1)}-g_n^{(1)}\|_{n}^2 \| \Phi_n \tilde{g}_m^{(2)}-g_n^{(2)}\|_{n}^2 \nn \\ 
&\leq 2 \, \| \Phi_n \tilde{g}_m^{(1)}-g_n^{(1)}\|_{n}^2 \left( \| \Phi_n \tilde{g}_m^{(2)}\|_{n}^2 + \| g_n^{(2)}\|_{n}^2 \right) \nn \\
&= 2 \, \| \Phi_n \tilde{g}_m^{(1)}-g_n^{(1)}\|_{n}^2 \left( \| \Phi_n \tilde{g}_m^{(2)}\|_{n}^2 + \frac{1}{\norm{f_n^{(2)}}_n^2 } \right) \nn 
\end{align}
and conclude from \eqref{Strong2condgm} and the convergence $\| \Phi_n \tilde{g}_m^{(2)}\|_{n} \to \|  \tilde{g}_m^{(2)}\|_{\infty}$ and the fact that $f^{(2)} \neq 0$.

To conclude that $\left\{ \,  f_n^{(1)}  \mid n \in \N \, \right\}$  converges weakly to $f^{(1)}$ it is enough to see that for $n$ large enough
\begin{align}
\left \langle f_n^{(1)}, g_n^{(1)} \right \rangle_{n,1} =   \left \langle (f_n^{(1)}, f_n^{(2)}) , (g_n^{(1)}, g_n^{(2)}) \right \rangle_{n,2}.
\end{align}
\epr 
Now we can show condition Mosco I. Let $\left\{ \,  \bif_n=\otimes_{j=1}^k f_n^{(j)} \mid n \in \N \, \right\}$ be a sequence of elements of the Hilbert space $L^2(\Omega_n^{(k)}, \mu_n^{(k)})$ converging weakly to a non-zero element $\bif=\otimes_{j=1}^k f^{(j)}$ of $ L^2(\Omega^{(k)}, \mu^{(k)})$. By Proposition~\ref{PropweakConv1fromk} we know that every $f_n^{(j)} $ converges weakly to $f^{(j)}$. Moreover, each function of the form $\mathbf{f}_n^{(-j)}=(f_n^{(1)}, \ldots, f_n^{(j-1)}, f_n^{(j+1)}, \ldots, f_n^{(k)})$ converges weakly towards  $\mathbf{f}^{(-j)}=(f^{(1)}, \ldots, f^{(j-1)}, f^{(j+1)}, \ldots, f^{(k)})$. Hence
\begin{align}
\liminf_{n \to \infty} \caE_n^{(k)}(\bif_n) &= \liminf_{n \to \infty} \sum_{j=1}^k \norm{\otimes_{\substack{l=1 \\ l \neq j}}^k  f_n^{(l)}}_{n,k-1}^2 \caE_n(f_n{(j)}) \nn \\
&\geq  \sum_{j=1}^k \norm{\otimes_{\substack{l=1 \\ l \neq j}}^k  f^{(l)}}_{\infty,k-1}^2 \caE(f^{(j)}) = \caE^{(k)}(\bif) .
\end{align}

\subsection{Proof of condition Mosco II}

\bpr
By Theorem 3.5 of \cite{andres2010particle}, Theorem 2.4 of \cite{kuwae2003convergence} and part 3 of Assumption~\ref{MainAssumption}, it is enough to verify \eqref{mosco2} for elements of $\Cfrak^{(k)}$. Let $\bif$ be such an element, i.e., there exists $M \in \N$ and $\left\{ \, c_i \in \R \mid i \in [k] \, \right\}$ such that 
\be 
\bif = \sum_{i=1}^{M} c_i \otimes_{j=1}^k f_{(i,j)}, 
\ee 
for some $f_{(i,j)} \in \Cfrak$. Moreover, by part 2 of Assumption~\ref{MainAssumption} we now that for each $f_{(i,j)} \in \Cfrak$, there exists a sequence $\left\{ \,  f_{(i,j)}^{(n)}\mid n \in \N \, \right\}$ of elements of $L^2(\Omega_n, \mu_n)$ converging strongly to $f_{(i,j)}$ and such that
\be\label{MoscoIfijn} 
\caE (f_{(i,j)}) = \lim_{n \to \infty} \caE_n (  f_{(i,j)}^{(n)} ).
\ee 
Now define the functions $\bif_n \in L^2(\Omega_n^{(k)}, \mu_n^{(k)})$ by
\be 
\bif_n := \sum_{i=1}^{M} c_i \otimes_{j=1}^k f_{(i,j)}^{(n)}. 
\ee 
We now show that the sequence $\left\{ \,  \bif_n \mid n \in \N \, \right\}$ converges strongly to $\bif$. From the strong convergence of $\left\{ \,  f_{(i,j)}^{(n)}\mid n \in \N \, \right\}$ to $f_{(i,j)}$ we know that there exists a sequence $\left\{ \,  \tilde{f}_{(i,j)}^{(m)}\mid m \in \N \, \right\}$ of elements of $\Cfrak$ such that
\be\label{Strong1cond}
\lim_{m \to \infty} \norm{ \tilde{f}_{(i,j)}^{(m)}-f_{(i,j)}}_{\infty} = 0 
\ee
and
\be\label{Strong2cond}
\lim_{m \to \infty} \limsup_{n \to \infty} \| \Phi_n \tilde{f}_{(i,j)}^{(m)}-f_{(i,j)}^n \|_{n} = 0. 
\ee
This suggests to consider the elements $\tilde{\bif}_m \in \Cfrak^{(k)}$ given by
\be 
\tilde{\bif}_m:= \sum_{i=1}^{M} c_i \otimes_{j=1}^k \tilde{f}_{(i,j)}^{(m)}. 
\ee 
To show the strong convergence  $\left\{ \,  \bif_n \mid n \in \N \, \right\} \to \bif$ it is enough to show that $\tilde{\bif}_m$ satisfies \eqref{DefStrongcond1} and \eqref{DefStrongcond2}. We only show \eqref{DefStrongcond1} since \eqref{DefStrongcond2} is done analogously. We estimate as follows.
\begin{align}
    \norm{\tilde{\bif}_m-\bif}_{\infty,k}^2 &= \norm{\sum_{i=1}^{M} c_i \otimes_{j=1}^k \tilde{f}_{(i,j)}^{(m)}-\sum_{i=1}^{M} c_i \otimes_{j=1}^k \tilde{f}_{(i,j)}}_{\infty,k}^2 = \norm{\otimes_{j=1}^k \sum_{i=1}^{M} c_i  \left(\tilde{f}_{(i,j)}^{(m)}- \tilde{f}_{(i,j)}\right)}_{\infty,k}^2 \nn \\
    &= \prod_{j=1}^k \norm{\sum_{i=1}^{M} c_i  \left(\tilde{f}_{(i,j)}^{(m)}- \tilde{f}_{(i,j)}\right)}_{\infty}^2  \leq \prod_{j=1}^k \sum_{i=1}^{M} |c_i|^2 \, \norm{ \tilde{f}_{(i,j)}^{(m)}- \tilde{f}_{(i,j)}}_{\infty}^2.
\end{align}
\br\label{RemarkLowerStrongConv}
Notice that the same procedure also shows the strong convergence of functions like $\otimes_{\substack{l=1 \\ l \neq j}}^k \sum_{i=1}^{M} c_i ^k f_{(i,l)}^{(n)}$ towards $\otimes_{\substack{l=1 \\ l \neq j}}^k \sum_{i=1}^{M} c_i ^k f_{(i,l)}$
\er
By \eqref{Strong1cond} we obtain \eqref{DefStrongcond1}. This together with an analogous estimate for \eqref{DefStrongcond2} shows that indeed $\left\{ \,  \bif_n \mid n \in \N \, \right\} \to \bif$ strongly. Notice
\be
\caE_n^{(k)}(\bif_n) = \sum_{j=1}^k \norm{\otimes_{\substack{l=1 \\ l \neq j}}^k \sum_{i=1}^{M} c_i ^k f_{(i,l)}^{(n)}}_{n,k-1}^2 \caE_n(f_{(i,j)}^{(n)}),
\ee 
and conclude the proof of Mosco II by \eqref{MoscoIfijn} and Remark~\ref{RemarkLowerStrongConv}.
\epr 

\section{Applications}

\subsection{Mosco convergence of $k$ independent random walkers}
Consider a single particle performing a simple symmetric random walk in the rescaled integer lattice $\Z_n:= \frac{1}{n}\Z$. This a Markov process with infinitesimal generator given by
\be\label{LdifftIRW2}
L^{\text{irw}}_n f(x) = \alpha \Delta_n f (x) =\frac{\alpha n^2}{2}\(f(x+ \tfrac{1}{N}) - 2 f(v) + f(x-\tfrac{1}{N})\), \qquad  x \in   \Z_n.
\ee
The sequence of Dirichlet forms $\lbrace \caE_n^{\text{rw}} : n \in \N   \rbrace$ associated to this processes is given by
\be\label{DirIRW}
\caE_n^{\text{irw}} (f) = - \alpha \sum_{x \in \Z_n} f(x) \Delta_n f (x) \, \mu_n(x).
\ee
where $\mu_n$ is $\frac{1}{n}$ times the counting measure. 
We know, see Appendix in \cite{ayala2021condensation}, that sequence of Dirichlet forms $\{(\caE_n^{\text{irw}},D(\caE_n^{\text{rw}}))\}_n$ converges to  the Brownian Dirichlet form $(\caE^{\text{bm}},D(\caE^{\text{bm}}))$, i.e. the Dirichlet form associated to the Brownian motion in $\R$
\be
\caE^{\text{bm}}(f)= \frac{\alpha}{2} \int_{\R} f\myprime(x)^2 dx.
\ee
where $D(\caE^{\text{bm}})$ is the Sobolev space of $\R$ of order 1.

This convergence is valid under the convergence of Hilbert spaces given in terms of the sequence of Hilbert spaces
\be\label{HrwN}
H^{\text{irw}}_n :=L^2 (\Z_n, \mu_n ),
\ee
and the Hilbert space
\be\label{Hrw}
H^{\text{bm}}:=L^2 (\R,dx)
\ee
 i.e. the space of Lebesgue square-integrable functions in $\R$, by means  of the  restriction operators
\be\label{PHINIRW}
\{\Phi_n:C_k^{\infty}(\R) \subset H^{\text{bm}} \to H^{\text{rw}}_n\}_n \qquad \text{defined by}\qquad    \Phi_n f = f \mid_{\Z_n}.
\ee

\br
The choice of the space of all compactly supported smooth functions 
\[
C:=C_k^{\infty}(\R)
\]
as dense set for our Hilbert space turns out to be particularly convenient since it is a core  of the Dirichlet form associated to the standard Brownian motion.
\er

A straight forward application of Theorem~\ref{MainThmkpart} in this context is the Mosco convergence of the Dirichlet forms associated to the $k$-particle independent motion
\begin{align}\label{DefIRWkDir}
\caE_n^{\text{irw},(k)}(f) &:= \frac{\alpha n^2 }{2 n^{k}} \sum_{i=1}^k \sum_{\bix \in \Z_n^{k}} \sum_{\sigma = \pm 1} \left( f(\bix^{i,\sigma/n}) -f(\bix)  \right)^2, 
\end{align}
to the Dirichlet form corresponding to $k$ independent Brownian particles
\be\label{DefBMkDir}
\caE_n^{\text{bm},(k)}(f) = \alpha \sum_{i=1}^k  \int_{\R^k}\left( \frac{\partial}{\partial x_i}f(\bix)\right)^2   \, d\bix 
\ee
with domain $D(\caE^{\text{bm},(k)})$ being the Sobolev space of $\R^k$ of order 1.
\begin{corollary}
For every $k \in \N$, under the Hilbert convergence of Lemma~\ref{LemmaConvKHilb}, the sequence of Dirichlet forms $\{(\caE_n^{\text{irw},(k)},D(\caE_n^{\text{rw},(k)}))\}_n$ given by \eqref{DefIRWkDir} converges, in the  sense of Mosco convergence of Dirichlet forms, to the  Brownian Dirichlet form $(\caE^{\text{bm},(k)},D(\caE^{\text{bm},(k)}))$ given by \eqref{DefBMkDir}.
\end{corollary}

\subsection{Symmetric Inclusion Process}
The Symmetric Inclusion Process of parameter $\alpha \in \R_+$ ($\SIP(\alpha)$) is an interacting particle system where particles randomly hop on in a Markovian way on the lattice $\Z_n:= \frac{1}{n} \Z$. Configurations of particles are denoted by the Greek letters $\eta$ and $\xi$, and are elements of $\Omega_n := \N^{\mathbb{Z}_n}$. This means that for all site $x \in \Z_n$ the variable $\eta(x)$ denotes the number of particles at that site. We describe the time evolution of the processes $\lbrace \eta_t : t \geq 0 \rbrace$ in terms of the following infinitesimal generator:
\begin{align}\label{SIPgen}
\loc_n f(\eta) &= n^{2}  \sum_{x \in \Z} \, \eta(x/n) \left( \alpha + \eta((x+1)/n) \right) \left(  f( \eta^{x/n,(x+1)/n}) - f(\eta) \right) \nn \\
&+n^{2} \sum_{x \in \Z} \, \eta(x/n) \left( \alpha + \eta((x-1)/n) \right) \left(  f( \eta^{x/n,(x-1)/n}) - f(\eta) \right).
\end{align}
This process is known to satisfy a self-duality relation which allows to compute moments of the time-evolved random variables $\eta_t(x)$ in terms of a simpler version of the process $\lbrace \eta_t : t \geq 0 \rbrace$ considered in coordinate notation. To be more precise, let us introduce the function $D \colon \Z_n^k \times \Omega_n  \to \R$ given by
\begin{equation}
\label{eq:dual-d}
D(\bix ,\eta)=\prod_{y\in\Z_n} d(\xi(\bix)(y),\eta(y)),
\end{equation}
where the $k$-particle configuration $\xi(\bix)$ is given by
\be
\xi(\bix) :=\(\xi(\bix)(y), y\in \Z_n\)\qquad \text{with} \qquad \xi(\mathbf x)(y)=\sum_{i=1}^k \1_{x_i=y},
\ee
and single-site self-duality functions are given by
\be\label{classicalsingleSIP}
d(m,n) := \1_{\{m \leq n\}} \, \frac{n!}{(n-m)!} \frac{\Gamma(\alpha)}{\Gamma(\alpha+m)}.
\ee
The self-duality relation, at coordinate level, reads 
\be\label{dual1coord}
\E_\eta \big[D(\xi,\eta_t)\big]= \E_{\mathbf x} \big[D(\xi(X^{(k)}(t)), \eta)\big],
\ee
where $\lbrace X^{(k)}(t): t \geq 0 \rbrace$ is the process in $\Z_n^k$ with generator
\be\label{SIPkgencoord}
L_n^{\text{sip},(k)} f(\mathbf x) = n^2 \sum_{i=1}^{k} \sum_{\sigma = \pm 1} \Bigg( \alpha +  \sum_{\substack{j=1\\ j \neq i}}^k \mathbf 1_{x_j=x_i +\sigma/n} \Bigg) \( f(\mathbf x^{i,i+\sigma}) -f(\mathbf x) \),
\ee
where $\mathbf x^{i,i+\sigma}$ denotes $\mathbf x$ after moving the particle in position $x_i$ to position $x_i+\sigma/n \in \Znd$.

\subsection{Mosco convergence of $k$-\SIP particles}
Let us first start with the following observation

\bp\label{revermpos}
The process $\lbrace X^{(k)}(t): t \geq 0 \rbrace$ is reversible with respect to the probability measure
\be\label{revermusip}
\nu(\mathbf{x}) =  \prod_{j \in \Z } \frac{\Gamma( \alpha + \sum_{i=1}^{k} \1_{x_i = j} )}{\Gamma(\alpha) }.
\ee
\ep

\bpr
By detailed balance it is enough to verify that $\mu$ satisfies the relation:
\be
\nu (\mathbf{x}) \( \alpha + \sum_{ \substack{ j=1 \\ j \neq i}}^{k} \1_{x^i +r = x^j} \) =  \nu (\mathbf{x}^{i,i+r}) \( \alpha + \sum_{ \substack{ j=1 \\ j \neq i}}^{k} \1_{x^i = x^j} \),
\ee
which is a consequence of the basic property of the Gamma function
$\Gamma (z+1) = z \Gamma(z)$.
\epr

As a consequence of Proposition~\ref{revermpos} we have that, for every $k \in \N$, the sequence of Dirichlet forms $\{(\caE_n^{\text{irw},(k)},D(\caE_n^{\text{rw},(k)}))\}_n$ corresponding to $k-\SIP$ particles is given by
\begin{align}\label{DefSIPkDir}
\caE_n^{\text{sip},(k)}(f) &:= \frac{\alpha n^2 }{2 n^{k}} \sum_{i=1}^k \sum_{\bix \in \Z_n^{k}} \sum_{\sigma = \pm 1} \left( f(\bix^{i,\sigma/n}) -f(\bix)  \right)^2 \, \nu(\bix) \nn \\
&+\frac{n^2 }{2 n^{k}} \sum_{i=1}^k \sum_{\sigma = \pm 1} \sum_{\bix \in \Z_n^{k}} \sum_{\substack{j=1 \\ i \neq j}}^k \mathbf 1_{x_j=x_i +\sigma/n}  \left( f(\bix^{i,\sigma/n}) -f(\bix)  \right)^2 \, \nu(\bix)
\end{align}
We then have the following iterated consequence of Theorem~\ref{MainThmkpart}.

\begin{corollary}\label{SIPcoroll}
For every $k \in \N$, under the Hilbert convergence of Lemma~\ref{LemmaConvKHilb}, the sequence of Dirichlet forms $\{(\caE_n^{\text{sip},(k)},D(\caE_n^{\text{sip},(k)}))\}_n$ given by \eqref{DefIRWkDir} converges, in the  sense of Mosco convergence of Dirichlet forms, to the  Brownian Dirichlet form $(\caE^{\text{bm},(k)},D(\caE^{\text{bm},(k)}))$ given by \eqref{DefBMkDir}.
\end{corollary}

\bpr
For simplicity of exposition let us assume that $\alpha \geq 1$, it is then clear that
\be\label{ComparingDirForms}
\caE_n^{\text{sip},(k)}(f) \geq \caE_n^{\text{irw},(k)}(f),
\ee
for all $f \in H_n^{(k)}$. This relation, immediately implies that the condition Mosco I holds for the Dirichlet forms $\{(\caE_n^{\text{sip},(k)},D(\caE_n^{\text{sip},(k)}))\}_n$ and $(\caE^{\text{bm},(k)},D(\caE^{\text{bm},(k)}))$. Moreover, for any $f \in C^\infty_c(\R^k)$, since the difference of the two forms has only contributions from lower dimensional hyper-planes (i.e. places where coordinates of $\bix$ coincide), it is also clear that
\be
\lim_{n \to\infty} \left| \caE_n^{\text{sip},(k)}(\Phi_n^{(k)} f) - \caE_n^{\text{irw},(k)}(\Phi_n^{(k)}f) \right| = 0,
\ee 
which imply condition Mosco II for the Dirichlet forms $\{(\caE_n^{\text{sip},(k)},D(\caE_n^{\text{sip},(k)}))\}_n$ since the functions $\Phi_n^{(k)}f$ converge strongly to $f \in C^\infty_c(\R^k)$ by Remark~\ref{StrongConvPhin}.
\epr 

\subsection{$k$-point correlations}
In this application we go back to the symmetric inclusion process in configuration notation. Let us start with a slight restriction on the initial configuration of particles of the process.
\begin{assumption}\label{HydroAssump}
The process $\{\,  \eta_{r(n)t} \, \}_{t \geq 0}$ is such that for all $n \in \N$ the initial measure $\nu_{\rho}^n$ is associated to a profile $\rho \in S(\R)$ and satisfies
\begin{itemize}
    \item There exists constants $C_k >0$,  independent of $n$, and $\epsilon >0$ such that 
    \be 
    \sup_{\bix \in \Z_n^k} \left| \int   D\left(\xi(\bix), \eta \right) \, \nu_{\rho}^n(d\eta) - \prod_{i=1}^k \int  \eta(x_i) \,  \nu_{\rho}^n(d\eta) \right| \leq C_k n^{-\epsilon}
    \ee 
    for all $k \in \N$.
\end{itemize}
\end{assumption}


We then have the following consequence of Corollary~\ref{SIPcoroll}.

\begin{corollary}\label{MoscoKpointsCor}
Under Assumption \ref{HydroAssump} we have that for all $t \geq 0$, and for all  $F \in S(\R^k)$
\be 
\lim_{n \to \infty} n^{-k} \, \sum_{\bix \in \Z_n^k}  F(\bix) \, \E_n\left[  D\left(\xi(\bix), \eta(n^2t) \right)\right] = \int_{\R^k} F(\biz) \, S_t^{\text{bm},(k)} \prod_{i=1}^k \rho(z_i) \, d\biz, 
\ee 
where $S_t^{\text{bm},(k)}$ denotes the semigroup of $k$ independent Brownian motions.
\end{corollary}

\bpr 
By self-duality we have that for all $\bix\in \Z_n^k$
\begin{align}
 \E_n\left[  D\left(\xi(\bix), \eta(n^2t) \right) \right] &= \int_\Omega \E_\eta \left[  D\left(\xi(\bix), \eta(n^2t) \right) \right]    \, \nu_\rho(d \eta)   \nn \\ 
 &= \int_\Omega \E_\bix  \left[  D\left(\xi(X^{(k)}(t)), \eta \right) \right]    \, \nu_\rho(d \eta)   \nn \\ 
  &= \sum_{\biy \in \Z_n^k} p_{n^2 t}^{\text{sip},(k)} (\bix,\biy) \, \int_\Omega  D\left(\xi(\biy, \eta \right)     \, \nu_\rho(d \eta),  
\end{align}
where $p_{n^2 t}^{\text{sip},(k)} (\bix,\biy)$ denotes the transition probability the process with generator \eqref{SIPkgencoord}. Then we have
\begin{align}
 &\left|   n^{-k} \, \sum_{\bix \in \Z_n^k}  F(\bix) \,  \E_n\left[  D\left(\xi(\bix), \eta(n^2t) \right)\right] -\int_{\R^k} F(\biz) S_t^{\text{bm},(k)} \prod_{i=1}^k \rho(z_i) \, d\biz \right| \nn \\
 &=\left|   n^{-k} \, \sum_{\bix \in \Z_n^k}  F(\bix) \,  \sum_{\biy \in \Z_n^k} p_{n^2 t}^{\text{sip},(k)} (\bix,\biy) \, \int_\Omega  D\left(\xi(\biy, \eta \right)     \, \nu_\rho(d \eta)  -\int_{\R^k} F(\biz) S_t^{\text{bm},(k)} \prod_{i=1}^k \rho(z_i) \, d\biz \right| \nn \\
&\leq  n^{-k} \, \sum_{\bix \in \Z_n^k}  F(\bix) \,  \sum_{\biy \in \Z_n^k} p_{n^2 t}^{\text{sip},(k)} (\bix,\biy) \,  \left| \int_\Omega  D\left(\xi(\biy, \eta \right)     \, \nu_\rho(d \eta)  - \prod_{i=1}^k \int  \eta(y_i) \,  \nu_{\rho}^n(d\eta)\right| \nn \\
&+\left|   n^{-k} \, \sum_{\bix \in \Z_n^k}  F(\bix) \,  \sum_{\biy \in \Z_n^k} p_{n^2 t}^{\text{sip},(k)} (\bix,\biy)\, \prod_{i=1}^k \int  \eta(y_i) \,  \nu_{\rho}^n(d\eta)     \, \nu_\rho(d \eta)  -\int_{\R^k} F(\biz) S_t^{\text{bm},(k)} \prod_{i=1}^k \rho(z_i) \, d\biz \right| \nn \\
&\leq  \frac{C_k}{n^{k+\epsilon}}\, \sum_{\bix \in \Z_n^k}  F(\bix) \, +\left|   n^{-k} \, \sum_{\bix \in \Z_n^k}  F(\bix) \,  \sum_{\biy \in \Z_n^k} p_{n^2 t} (\bix,\biy)\, \prod_{i=1}^k \rho(y_i)  -\int_{\R^k} F(\biz) S_t^{\text{bm},(k)} \prod_{i=1}^k \rho(z_i) \, d\biz \right|,
\end{align}
and conclude by Corollary~\ref{SIPcoroll} and Theorem~\ref{MKS}.
\epr 

\br
Notice that in the presence of weak convergence of the $k$-particle coordinate $\SIP$ towards $k$ independent Brownian particles we have
\be 
\lim_{n \to \infty}  \E_n\left[  D\left(\xi(\bix^{(n)}), \eta(n^2t) \right)\right] =  S_t^{\text{bm},(k)} \prod_{i=1}^k \rho(x_i),  
\ee 
where for any $\bix \in \R$ we used the notation $\bix^{(n)}:= \tfrac{1}{n}\lfloor n \bix \rfloor \in \Z_n^k$. As consequence Corollary~\ref{MoscoKpointsCor} can also be derived. However to derive Corollary~\ref{MoscoKpointsCor} it is enough to have the strong converge 
\be 
\E_n\left[  D\left(\xi(\bix^{(n)}), \eta(n^2t) \right)\right] \to   S_t^{\text{bm},(k)} \prod_{i=1}^k \rho(x_i)  
\ee 
in the sense of Definition \ref{strongcon}.
\er 

\section*{Acknowledgements}
The author expresses sincere gratitude to J. Zimmer for their invaluable discussions.

\bibliography{biblio} 
\bibliographystyle{plain} 
\end{document}

%% file: mymacros.tex
\newcommand{\Z}{\mathbb Z}
\newcommand{\R}{\mathbb R}
\newcommand{\N}{\mathbb N}

\newcommand{\E}{\mathbb E}

\newcommand{\Znd}{\mathbb Z_n^d}

\renewcommand{\phi}{\varphi}

\newcommand{\loc}{\mathcal{L}}

\newcommand{\Cfrak}{{\mathfrak C}}

\newcommand{\biy}{{\mathbf y}}
\newcommand{\biz}{{\mathbf z}}
\newcommand{\bif}{{\mathbf f}}
\newcommand{\bfg}{{\mathbf g}}

\def\1{{\mathchoice {\rm 1\mskip-4mu l} {\rm 1\mskip-4mu l}
{\rm 1\mskip-4.5mu l} {\rm 1\mskip-5mu l}}}

\newtheorem{theorem}{{\small T}{\scriptsize HEOREM}}[section]
\newtheorem*{theorem*}{{\small T}{\scriptsize HEOREM}}
\newtheorem{corollary}{{\bf{\small C}{\scriptsize OROLLARY}}}[section]
\newtheorem{proposition}{{\bf{\small P}{\scriptsize ROPOSITION}}}[section]
\newtheorem*{proposition*}{{\bf{\small P}{\scriptsize ROPOSITION}}}
\newtheorem{lemma}{{\bf{\small L}{\scriptsize EMMA}}}[section]
\newtheorem{remark}{{\bf{\small R}{\scriptsize EMARK}}}[section]
\newtheorem{definition}{{\bf{\small D}{\scriptsize EFINITION}}}[section]
\newtheorem{induction}{{\bf{\small I}{\scriptsize NDUCTIVE HYPOTHESIS}}}[section]
\newtheorem{assumption}{{\bf{\small A}{\scriptsize SSUMPTION}}}[section]

\renewenvironment{proof}[1]
{\noindent{{\bf{\small{ P}{\scriptsize ROOF}}}.}\hspace{0.1cm} #1} {$\;\qed$\newline}

\newcommand{\beq}{\begin{eqnarray}}
\newcommand{\eeq}{\end{eqnarray}}

\newcommand{\ba}{\begin{align*}}
\newcommand{\ea}{\end{align*}}

\newcommand{\be}{\begin{equation}}
\newcommand{\ee}{\end{equation}}

\newcommand{\bl}{\begin{lemma}}
\newcommand{\el}{\end{lemma}}

\newcommand{\br}{\begin{remark}}
\newcommand{\er}{\end{remark}}

\newcommand{\bt}{\begin{theorem}}
\newcommand{\et}{\end{theorem}}

\newcommand{\bd}{\begin{definition}}
\newcommand{\ed}{\end{definition}}

\newcommand{\bind}{\begin{induction}}
\newcommand{\eind}{\end{induction}}

\newcommand{\basu}{\begin{assumption}}
\newcommand{\esu}{\end{assumption}}

\newcommand{\bp}{\begin{proposition}}
\newcommand{\ep}{\end{proposition}}

\newcommand{\bc}{\begin{corollary}}
\newcommand{\ec}{\end{corollary}}

\newcommand{\bpr}{\begin{proof}}
\newcommand{\epr}{\end{proof}}

\newcommand{\bi}{\begin{itemize}}
\newcommand{\ei}{\end{itemize}}

\newcommand{\ben}{\begin{enumerate}}
\newcommand{\een}{\end{enumerate}}

\newcommand{\caE}{{\mathrsfs E}}

\newcommand{\SIP}{\text{\normalfont SIP}}

\newcommand{\bix}{{\mathbf x}}

\renewcommand{\(}{\left(}
\renewcommand{\)}{\right)}
\newcommand{\nn}{\nonumber}

\newmuskip\pFqmuskip

\newcommand\pFq[6][8]{%
	\begingroup 
	\pFqmuskip=#1mu\relax
	\mathcode`\,=\string"8000
	\begingroup\lccode`\~=`\,
	\lowercase{\endgroup\let~}\pFqcomma
	{}_{#2}F_{#3}{\left[\genfrac..{0pt}{}{#4}{#5};#6\right]}%
	\endgroup
}
\newcommand{\pFqcomma}{\mskip\pFqmuskip}

\newcommand*{\myprime}{^{\prime}\mkern-1.2mu}

\newcommand{\norm}[1]{\left\lVert#1\right\rVert}

\newcommand{\xangle}{11}
\newcommand{\yangle}{133}
\newcommand{\zangle}{270}

\newcommand{\xlength}{1}
\newcommand{\ylength}{1}
\newcommand{\zlength}{1}

\pgfmathsetmacro{\xx}{\xlength*cos(\xangle)}
\pgfmathsetmacro{\xy}{\xlength*sin(\xangle)}
\pgfmathsetmacro{\yx}{\ylength*cos(\yangle)}
\pgfmathsetmacro{\yy}{\ylength*sin(\yangle)}
\pgfmathsetmacro{\zx}{\zlength*cos(\zangle)}
\pgfmathsetmacro{\zy}{\zlength*sin(\zangle)}